\documentclass[12pt]{amsart}
\usepackage{amssymb,amsmath,amsthm}
\usepackage{a4}
\usepackage{color}
\usepackage[colorlinks=true,
linkcolor=webgreen,
filecolor=webbrown,
citecolor=webgreen]{hyperref}
\definecolor{webgreen}{RGB}{0,0,1}
\definecolor{recrown}{RGB}{1,.2,.6}

\begin{document}
\newtheorem{theorem}{Theorem}
\newtheorem{corollary}[theorem]{Corollary}
\newtheorem{lemma}[theorem]{Lemma}
\theoremstyle{definition}
\newtheorem{example}{Example}
\newtheorem*{examples}{Examples}
\newtheorem*{notation}{Notation}
\newtheorem*{remark}{Remark}
\theoremstyle{theorem}
\newtheorem{thmx}{Theorem}
\renewcommand{\thethmx}{\text{\Alph{thmx}}}
\newtheorem{lemmax}{Lemma}
\renewcommand{\thelemmax}{\text{\Alph{lemmax}}}
\leftmargin=.5in
\rightmargin=0.5in
\theoremstyle{definition}
\newtheorem*{definition}{Definition}
\title[]{\bf Another irreducibility criterion}
\author{Jitender Singh$^{1}$}
\address[1]{Department of Mathematics,
Guru Nanak Dev University, Amritsar-143005, India\linebreak
 {\tt sonumaths@gmail.com}}
\author{Sanjeev Kumar$^{2,*}$}
\address[2]{Department of Mathematics,
SGGS College, Sector-26, Chandigarh-160019, India\linebreak
{\tt sanjeev\_kumar\_19@yahoo.co.in}}
\markright{}
\date{}
\maketitle
\parindent=0cm
\footnotetext[1]{Corresponding Author email: {\tt sanjeev\_kumar\_19@yahoo.co.in}\\

2010MSC: {Primary 12E05; 11C08}\\

}
\begin{abstract}
  Let $f=a_0+ a_{1}x+\cdots+a_m x^m\in \Bbb{Z}[x]$ be a primitive polynomial. Suppose that there exists a positive real number  $\alpha$ such that
$|a_m| \alpha^m>|a_0|+|a_1|\alpha+\cdots+|a_{m-1}|\alpha^{m-1}$. We prove that if there exist natural numbers $n$ and $d$ satisfying $n\geq \alpha+ d$ for which either $|f(n)|/d$ is a prime, or $|f(n)|/d$ is a prime-power coprime to $|f'(n)|$, then $f$ is irreducible in $\mathbb{Z}[x]$.
\end{abstract}
\section{Introduction.}
The classical irreducibility criteria  due to Sch\"onemann (1846), Eisenstein (1850), Dumas (1906), and Perron (1907) have become paradigm for testing irreducibility of polynomials having rational coefficients. The demesne revealing riveting facts about irreducibility of polynomials over
prescribed domains has always been the cradle of such baroque classical results which for decades have witnessed
cogent extensions and generalizations. Such irreducibility criteria have exhibited a close affinity
to prime numbers and primality as is evident from the illustrious Buniakowski's conjecture of 1854 which asserts that if $f$ is an irreducible polynomial having integer coefficients such that the elements in the set $f(\mathbb N)$ have no common factors other than $\pm 1$, then the set $f(\mathbb N)$ contains infinitely many  prime numbers. The converse of Buniakowski's conjecture holds affirmatively via primality.

Another classical irreducibility result due to A. Cohn \cite[p.~133]{P} states that if a prime number can be expressed in base 10 as
$\sum_{i=0}^m a_i10^i$ for some positive integer $m$, then the polynomial $\sum_{i=0}^m a_ix^i$ is irreducible in $\mathbb{Z}[x]$. Cohn's result was then generalized to arbitrary base by Brillhart et al. \cite{B} and further in  Bonciocat et al. \cite{Bonciocat}. In \cite{Mu}, Murty provided elementary proof of Cohn's irreducibility criterion. Interestingly, one of the main results of Murty \cite{Mu}, generalized by Girstmair \cite{G}  apprised of a strong converse of Buniakowski's conjecture which was further generalized in \cite{Jakhar} and \cite{JSSK2022} for polynomials having integer coefficients.
\begin{thmx}[\cite{Jakhar}]\label{th:A}
    Let $f=a_0+ a_{1}x+\cdots+a_m x^m\in \Bbb{Z}[x]$ be a primitive polynomial. Suppose there exists a positive real number  $\alpha$ such that
\begin{eqnarray*}
|a_m| \alpha^m>|a_0|+|a_1|\alpha+\cdots+|a_{m-1}|\alpha^{m-1}.
\end{eqnarray*}
If there exist natural numbers $n$ and $d$ satisfying $n\geq \alpha+ d$ for which $f(n)=\pm pd$ for  a prime $p$, then $f$ is irreducible in $\mathbb{Z}[x]$.
 \end{thmx}
\begin{thmx}[\cite{JSSK2022}]\label{th:B}
  Let $f=a_0+a_1 x+\cdots+a_mx^m\in \Bbb{Z}[x]$ be primitive, and let
  \begin{eqnarray*}
H=\max_{0\leq i\leq m-1}\{|a_i/a_m|\}.
  \end{eqnarray*}
  Let $f'(x)$ denote the formal derivative of $f(x)$ with respect to $x$. If there exist natural numbers $n$, $d$, $k$, and a prime $p\nmid d$ such that $n\geq 1+H+d$, $f(n)=\pm p^k d$,  and for $k>1$, also $p\nmid f'(n)$, then $f$ is irreducible in $\Bbb{Z}[x]$.
 \end{thmx}
In the present note, we generalize Theorem \ref{th:A} to the case when $|f(n)|/d$ is a prime-power with the mild condition of coprimality of $|f(n)|/d$ with $|f'(n)|$. More precisely, we have the following result.
\begin{theorem}\label{th:1}
    Let $f=a_0+ a_{1}x+\cdots+a_m x^m\in \Bbb{Z}[x]$ be a primitive polynomial. Suppose that there exists a positive real number $\alpha$ such that
\begin{eqnarray*}
|a_m| \alpha^m>|a_0|+|a_1|\alpha+\cdots+|a_{m-1}|\alpha^{m-1}.
\end{eqnarray*}
If there exist natural numbers $n$ and $d$ satisfying $n\geq \alpha+ d$ for which  $|f(n)|/d$ is prime, or $|f(n)|/d$ is a prime-power coprime to $|f'(n)|$, then $f$ is irreducible in $\mathbb{Z}[x]$.
 \end{theorem}
 \begin{example}
For $k\geq m+2\geq 4$ and $p\geq 1+d$, the polynomial
    \begin{eqnarray*}
X=-p+x\pm (p^{k-m}d)x^m
    \end{eqnarray*}
  satisfies the hypothesis of Theorem \ref{th:1} with $\alpha=1$, $a_0=-p$, $a_1=1$; $a_i=0$ for $i=2,3,\ldots,m-1$;  $a_m=\pm p^{k-m} d$; and  $n=p\geq 1+d=\alpha+d$, since we have $X(p)=\pm p^{k} d$; $X'(p)\equiv 1\mod p$ so that $\gcd(|X(p)|/d,|X'(p)|)=1$, and
  \begin{eqnarray*}
|a_m|\alpha^m=p^{k-m} d\geq p^2>p+1=\sum_{i=0}^{m-1}|a_i|\alpha^i.
    \end{eqnarray*}
By Theorem \ref{th:1}, the polynomial $X$ is irreducible in $\mathbb{Z}[x]$.
 \end{example}
 \begin{example}
Now consider the polynomial
    \begin{eqnarray*}
Y=(x-p)+(x-p)^2+\cdots+(x-p)^{m-1}\pm (p^{2k-1}d)x^m
    \end{eqnarray*}
for $k\geq m\geq 2$ and $p\geq 1+d$. Here, $a_i=\sum_{j=i}^{m-1} {j\choose i}(-p)^{j-i}$ for $i=0,1,\ldots, m-1$; $a_m=\pm p^{2k-1} d$, $\alpha=1$, and  $n=p\geq 1+d$. We find that $Y(p)=\pm p^{2k+m-1} d$, $Y'(p)\equiv 1\mod p$. These along with the fact that $p^2>1+p$ yield the following:
  \begin{eqnarray*}
|a_m|\alpha^m=\frac{p^{2k}d}{p}\geq \frac{(p^2)^m d}{p}>\frac{(1+p)^m}{p}>(1+p)\frac{(1+p)^{m-1}-1}{1+p-1}=\sum_{i=0}^{m-1}|a_i|\alpha^i.
    \end{eqnarray*}
Since $a_{m-1}=1$, it follows that $Y$ is a primitive polynomial. By Theorem \ref{th:1}, the polynomial  $Y$ is irreducible in $\mathbb{Z}[x]$.
 \end{example}
\section{Proof of Theorem 1.}
Let $|f(n)|/d=p^k$ for some prime $p$ and positive integer $k$.  If $|x|\geq \alpha$, then in view of the hypothesis, we have $|a_m| \alpha^m>\sum_{j=0}^{m-1}|a_j|\alpha^j$. Consequently, we have
\begin{eqnarray*}
|f(x)|\geq |x|^m\Bigl(|a_m|-\sum_{i=0}^{m-1}|a_i||x|^{-(m-i)}\Bigr)\geq \alpha^m\Bigl(|a_m|-\sum_{i=0}^{m-1}|a_i|\alpha^{-(m-i)}\Bigr)>0,
\end{eqnarray*}
which shows that each zero $\theta$ of $f$ satisfies $|\theta|<\alpha$.

Now assume on the contrary that $f(x)=f_1(x)f_2(x)$ for nonconstant polynomials $f_1$ and $f_2\in \mathbb{Z}[x]$. Since we have
\begin{eqnarray*}
\pm p^kd=f(n)=f_1(n)f_2(n),
\end{eqnarray*}
at least one of $|f_1(n)|$ and $|f_2(n)|$ is divisible by $p$. Assume that $p$ divides $|f_2(n)|$.
Firstly, let us suppose that $p$ does not divide $|f_1(n)|$. Then $p^k$ divides $|f_2(n)|$, and so, $|f_1(n)|$ must divide $d$ so that we have $|f_1(n)|\leq d$. If $\beta~(\neq0)$ is the leading coefficient of $f_1$, then
\begin{eqnarray*}
f_1(n)=\beta\prod_\theta (n-\theta),
\end{eqnarray*}
where the product runs over all zeros $\theta$ of $f_1$. Observe that each such $\theta$ satisfies $|\theta|<\alpha$. Since
\begin{eqnarray*}
|n-\theta|\geq n-|\theta|>n-\alpha\geq d,
\end{eqnarray*}
 we arrive at the following:
\begin{eqnarray*}
d\geq |f_1(n)|=|\beta|\prod_{\theta}|n-\theta| &>& |\beta|d^{\deg f_1}\geq |\beta|d\geq d,
\end{eqnarray*}
leading to a contradiction.

Now assume that $p$ divides $|f_1(n)|$. Since $p$ divides $|f_2(n)|$, we must have $k\geq 2$. Consequently,  $p$ divides $|{f_1}'(n)f_2(n)+f_1(n){f_2}'(n)|$, which  in view of the fact that
 \begin{eqnarray*}
  {f_1}'(n)f_2(n)+f_1(n){f_2}'(n)=f'(n),
 \end{eqnarray*}
shows that $p$ divides $|f'(n)|$. This contradicts the hypothesis. So, $f$ must be irreducible in $\mathbb{Z}[x]$. \qed

The following remark and examples serve well to make the present idea efficaciously
comprehensible rendering an advantage over the results already known in the domain.
\begin{remark} Note that Theorem \ref{th:A} is the special case  of Theorem \ref{th:1} with $k=1$. The significance of Theorem \ref{th:1} lies in the fact that whenever each one of Theorems \ref{th:A}, \ref{th:B}, \ref{th:1} is applicable, Theorems \ref{th:A}, \ref{th:B} may encounter a tedious factorization of integers. This is demonstrated in the following explicit examples.
\end{remark}
\begin{example}
 Consider the polynomial
\begin{eqnarray*}
Z=9-x+72 x^{18}.
\end{eqnarray*}
The smallest value of $n$ for which Theorem \ref{th:1} is applicable for $Z$ is $n=9$ with $\alpha=1$, $d=8$, and $Y(9)/8=3^{38}$, whereas the smallest value of $n$ for which  Theorems \ref{th:A} and \ref{th:B} are applicable is $n=28$ with  $d=13$ and
\begin{eqnarray*}
Z(28)/13=
619774506599223645785433953,
\end{eqnarray*}
which is an 18-digit prime number.
\end{example}
\begin{example}
Consider the following polynomials $Z_d$ as mentioned in \cite{JSSK2022}
\begin{equation*}
Z_d=p^k-x\pm (p^k d) x^m,~2\leq d\leq p^k-1,~k\geq 2,
\end{equation*}
where $k,m,d$ are positive integers and $p$ is a prime number. Here, $a_0=p^k$, $a_1=-1$, $a_i=0$ for $i=2,\ldots,m-1$, and $a_m=\pm p^k d$.
Taking $\alpha=1$ and $n=p^k$, we have
\begin{eqnarray*}
|a_m|\alpha^m=p^k d>p^k+1&=&\sum_{i=0}^{m-1}|a_i|\alpha^i;~n=p^k\geq 1+d=\alpha+d,\\
|Z_d(p^k)|/d&=&p^{k(1+m)};~{Z_d}'(p^k)\equiv -1\mod p,
\end{eqnarray*}
so that $|Z_d(p^k)|/d$ is coprime to $|{Z_d}'(p^k)|$. Thus  by Theorem \ref{th:1}, the polynomial $Z_d$ is irreducible in $\mathbb{Z}[x]$.

Here, for the aforementioned value of $n$ and $\alpha$, $Z_{p^k-1}$ is irreducible by Theorem \ref{th:1}, the irreducibility of which cannot be easily concluded from Theorem \ref{th:A} or  Theorem \ref{th:B}.
\end{example}
\bibliographystyle{amsplain}

\end{document}